\documentclass[a4paper,12pt]{amsart}

\usepackage{changebar}
\let\cbstart\relax
\let\cbend\relax

\usepackage[latin1]{inputenc}
\usepackage[T1]{fontenc}
\usepackage{url}
\usepackage{amsmath,amssymb}
\IfFileExists{dsfont}{\usepackage{dsfont}\let\mathbb\mathds}{}
\usepackage[english]{babel}
\usepackage[all,dvips]{xy}
\usepackage{fleches}

\makeatletter
  \begingroup
    \divide\count0 by 60
    \multiply\count0 by 60
    \advance\count1 by -\count0
    \xdef\dateandtime{%
      \the\year/\two@digits{\the\month}/\two@digits{\the\day}\ 
      \ifnum\count2<10 0\fi \the\count2:%
      \ifnum\count1<10 0\fi \the\count1
    }%
  \endgroup
\makeatother
\usepackage{fancyhdr}
  \fancypagestyle{firstpage}{%
    \fancyhf{}
    \fancyfoot[CO,CE]{\footnotesize\thepage}
    
    }

\makeatletter
\def\@secnumfont{}
\def\section{\@startsection{section}{1}%
  \z@{.7\linespacing\@plus\linespacing}{.5\linespacing}%
  {\normalfont\bfseries\large\centering}}
\makeatother

\usepackage{amsthm}
\theoremstyle{plain}
\newtheorem{theo}{Theorem}[section]
\newtheorem{lemma}[theo]{Lemma}
\newtheorem{prop}[theo]{Proposition}
\newtheorem{coro}[theo]{Corollary}
\theoremstyle{definition}
\newtheorem{rem}[theo]{Remark}

\newtheorem{defi}[theo]{Definition}

\def\H{\mathscr{H}}
\def\T{\mathscr{T}}
\def\CC{\mathbb{C}}
\let\phi\varphi
\let\epsilon\varepsilon
\def\SymetricGroup{\mathfrak{S}}

\DeclareMathOperator{\spec}{Spec}
\DeclareMathOperator{\im}{Im}
\let\hom\relax
\DeclareMathOperator{\hom}{Hom}
\DeclareMathOperator{\gal}{Gal}
\DeclareMathOperator{\sk}{Sk}
\DeclareMathOperator{\cosk}{Cosk}
\DeclareMathOperator{\ob}{Ob}
\DeclareMathOperator{\pr}{pr}
\DeclareMathOperator{\aut}{Aut}
\DeclareMathOperator{\LC}{LC}
\DeclareMathOperator{\LCF}{LCF}
\DeclareMathOperator{\SLC}{SLC}
\DeclareMathOperator{\SLCF}{SLCF}
\newcommand{\proobj}{\mathop{{\text{``}}\,{\limpro}{\text{''}}}}

\def\op{^{\text{\rm op}}}
\def\et{_{\text{\rm ét}}}
\def\ie{\emph{i.e.}}
\let\tens\otimes
\let\wh\widehat

\let\iso\simeq
\let\equiv\approx
\let\nothing\varnothing
\def\pro{\text{\rm pro-}}
\let\psd\rtimes
\let\inter\cap

\def\d#1{\textbf{\textup{#1}}}

\usepackage{mathrsfs}
\def\fx#1{\mathscr{#1}}

\let\f\mathcal

\def\cat#1{\mathfrak{#1}}
\def\dcat#1{\boldsymbol{\mathfrak{#1}}}
\IfFileExists{font-lucidablackletter.sty}{%
  \usepackage{font-lucidablackletter}
  \def\cat##1{\text{\rm\csname textlucidablackletter\endcsname{##1}}}
  \def\dcat##1{\text{\rm\csname textlucidablackletter\endcsname{\fontseries{b}\selectfont##1}}}
}{}

\makeatletter
  \let\OLD@times\times
  \def\Times{\@ifnextchar_\timesWSB\OLD@times}
  \def\timesWSB_#1{\mathbin{\mathop{\OLD@times}\limits_{#1}}}
  \let\OLD@otimes\otimes
  \def\Otimes{\@ifnextchar_\otimesWSB\OLD@otimes}
  \def\otimesWSB_#1{\mathbin{\mathop{\OLD@otimes}\limits_{#1}}}
  
\makeatother
\newcommand{\definefunction}[4]{%
  \ensuremath{
    \left\{
      \begin{array}{ccl} 
        {\displaystyle #1} & \lra & {\displaystyle #2} \\
        {\displaystyle #3} & \lmt & {\displaystyle #4} \\
      \end{array}
    \right.}}
\newcounter{theoenumcounter}
\def\theoenumlabel{\textup{(}\textit{\roman{theoenumcounter}}\textup{)}}
\newenvironment{theoenum}{%
  \begin{list}{\theoenumlabel}{%
      \usecounter{theoenumcounter}%
      \topsep=2pt\itemsep=0pt}}
  {\end{list}}
\newcommand{\UN}[4][r]{%
  \ar@/^1pc/[#1]^{#2}_*=<0.3pt>{}="HAUT"
  \ar@/_1pc/[#1]_{#3}^*=<0.3pt>{}="BAS"
  \save\POS "HAUT",*{
      \vrule height 1pt depth 1pt width 0pt
      \vrule height 0pt depth 0pt width 1pt
      }="HAUT",\restore
  \save\POS "BAS",*{
      \vrule height 1pt depth 1pt width 0pt
      \vrule height 0pt depth 0pt width 1pt
      }="BAS",\restore
  \ar @{=>} "HAUT";"BAS" ^{#4}
}
\def\XYMATRIX#1{\raisebox{.5\depth}{\xymatrix{#1}}}

\begin{document}
\title[The Fundamental Group of an Algebraic Stack]{The Fundamental Group \\ of an Algebraic Stack}
\author{Vincent Zoonekynd}
\address{Université Paris 7,
UFR de Mathématiques,
Équipe Topologie et géométrie algébriques,
Case 7012,
2 place Jussieu,
75251 Paris Cedex 05}
\email{zoonek@math.jussieu.fr}
\urladdr{http://www.math.jussieu.fr/\textasciitilde zoonek/}
\begin{abstract}
  We define the notion of fundamental group of an algebraic
  stack, prove a comparison theorem between the fundamental
  group of a stack over $\CC$ and that of the associated
  analytic orbifold, show that this notion coincides with
  that defined with simplicial schemes by \cite{oda} and
  finally prove some short exact sequences involving these
  groups.
\end{abstract}
\maketitle

\section{Pointless fundamental groupoid of a topos}

We shall first recall some results from \cite{leroy}.

\begin{defi}
  A \d{groupoid} is a category all of whose morphisms are
  isomorphisms. As a group is a groupoid with a single
  object, a groupoid may conversely be seen as a ``group
  with several objects''.
  
  The \d{classifying topos} of a groupoid $G$ is the
  category of presheaves on the opposite category: 
  $\cat B G = G\op\cat{Set}$; in particular, the classifying
  topos of a group $G$ is its category of (left) $G$-sets. 

  A \d{progroupoid} is a projective $1$-system 
  $(G_i)_{i\in I}$ of groupoids, \ie, a functor 
  $I \lra \cat{Groupoids}$ from a
  filtering category $I$;
  we shall often denote it as
  $\proobj G_i$.

  A \d{strict progroupoid} is a progroupoid 
  $\proobj G_i$ such that 
  \begin{theoenum}
  \item the functors $G_i \lra G_j$ induce bijections
    between the sets of objects $\ob G_i \lra \ob G_j$;
  \item for any objects $X$, $Y$, the functors 
    $G_i \lra G_j$
    induce onto maps
    $\hom_{G_i} (X,Y) \lra \hom_{G_j} (X,Y)$.
  \end{theoenum}
  All our progroupoids will be strict.

  The classifying topos of a progroupoid 
  $\proobj G_i$ is 
  $ \cat B \proobj G_i := 
  \dlimpro\nolimits^{\dcat{Topos}}
  \cat B G_i. $
  (The superscript ``$\dcat{Topos}$'' indicates the
  $2$-category in which the $2$-limit is taken.)
\end{defi}

\begin{defi}
  Let $\T$ be a topos. 
  An object $X$ is said to be \d{connected} if whenever
  $X = X_1 \amalg X_2$, 
  either $X_1$ or $X_2$ is a final object (often denoted $\nothing$).
  A \d{constant object} is a disjoint sum of copies of the
  final object $*$, \ie, it may be written
  $I_\T = \coprod_{i\in I} *$ for some set $I$. 
  The topos $\T$ is said to be connected if its final object is. 

  An object $F$ is said to be \d{trivialized by} an object
  $U$ if its \d{restriction}
  $F|_U = (F \times U \lra U)$ is a constant object of the
  topos $\T/U$. An object $F$ is said to be locally constant
  if the objects trivializing it cover the final object, 

  A \d{Galois object} in a topos $\T$ is a connected,
  locally constant object~$Y$ such that the morphism
  $Y \times (\aut Y)_\T \lra Y \times Y$ 
  is an isomorphism. 

  A category $\cat C$ is said to be \d{generated by} a set
  of objects $\cat K \subset \ob\cat C$ if 
  for any objects $x$, $y$, for any distinct morphisms $f, g
  : x \lra y$, there exists an object $z \in \cat K$ and a
  morphism $k: z \lra x$ such that $fh \neq gh$.

  A \d{Galois topos} is a topos generated by its Galois
  objects. 
\end{defi}

We shall see that Galois toposes are exactly the classifying
toposes of strict groupoids. 

\begin{defi}
  A topos is said to be \d{locally connected} if it is
  generated by its connected objects; this is equivalent to
  require that the ``constant object'' functor
  $\cat{Set} \lra \T$ have a left adjoint $\pi$, which we
  shall call the \d{connected components functor}.
\end{defi}
\begin{defi}
  A \d{projective 2-system} in a $2$-category $\dcat A$ 
  is a $2$-functor 
  $I \lra \dcat{A}$ from a small filtering $1$-category
  $I$. A projective $2$-system of toposes is said to be
  \d{strict} if the inverse image 
  functors of the topos morphisms involved are fully
  faithful.

  Let $\dcat A$ be a $2$-category and $\dcat B$ a (full)
  sub-$2$-category. A 1-morphism $f:A \lra B$ to an object of
  $\dcat B$ is said to be \d{2-universal} if
  for any morphism 
  $g : A \lra B'$ to an object of $\dcat B$, 
  there exists a morphism 
  $h : B \lra B'$ and a $2$-isomorphism 
  $\alpha : h f \lra g$ such that for any other morphism 
  $h' : B \lra B'$ and $2$-isomorphism $\alpha' : 
  h' f \lra g$, there exists a \emph{unique} $2$-isomorphism
  $\beta : h \lra h'$ such that 
  $$ \xymatrix{ hf \ar[rd]^\alpha \ar[d]_{\beta f} \\ 
    h'f \ar[r]_{\alpha'} & g.} $$
\end{defi}
\begin{theo}
  Let $\T$ be a locally connected topos.
  \begin{theoenum}
  \item The full subcategory $\SLC\T$ of disjoint sums of
    locally constant objects of $\T$ is a Galois topos and
    the inclusion functor $\SLC\T\lra\T$ is the inverse
    image functor of a topos morphism
    $\T \lra \SLC\T$ which is $2$-universal among morphisms
    to Galois toposes.
  \item The Galois toposes are the strict
    projective $2$-limits (in the $2$-category of toposes) 
    of classifying toposes of groupoids. 
  \end{theoenum}
\end{theo}
\begin{rem}
  We shall need the following points of the proof of this
  theorem. 
  \begin{theoenum}
  \item If $R$ is a covering sieve of the final object $*$
    of $\T$ (we shall write: $R \in J(*)$),
    denote $\LC(\T,R)$ the full subcategory of $\T$
    of the objects trivialized by the connected objects of
    $R$: it is a topos and the inclusion functor 
    $\LC(\T,R) \lra \T$ is the inverse image functor of a
    topos morphism.
  \item The toposes $\LC(\T,R)$ are classifying toposes of
    groupoids; more precisely, if $G$ is the groupoid of
    points of $\LC(\T,R)$, \ie, the groupoid of topos
    morphisms $\cat{Set}\lra\LC(\T,R)$ (and their $2$-isomorphisms), 
    then 
    $ \LC(\T,R) \equiv \cat B G. $
  \item The Galois topos $\SLC\T$ is the projective
    $2$-limit of these toposes: 
    $$ \SLC\T \equiv \dlimpro_{R \in J(*)}
    \!\!{}^{\dcat{Topos}}
    LC(\T,R). $$
  \end{theoenum}
\end{rem}

\section{Pointed fundamental progroupoid of a topos}

The following lemma
 gives an explicit description of the
classifying topos of a progroupoid. In particular, the
connected objects of $\cat B \proobj G_i$ are the (images of
the) connected objects of the various $\cat B G_i$. 

\begin{lemma}
  Let $G = \proobj G_i$ be a strict progroup and 
  $\cat S_0$ the category whose objects are sets~$X$ endowed
  with an action of one of the $G_i$ (hence of any $G_k$ for
  $k \lra i$), connected as presheaves (\ie, they have a
  single orbit) and whose morphisms are 
  $$
  \raisebox{.5\depth}{\xymatrix @!0 @R=.5em @C=3pc{
      & i \\
      *!!<0pt,0.7ex>+!R!<-1.2ex,0pt>{\hom\bigl( (X, G_i), (Y, G_j) \bigl) = 
        \hom_{G_k} (X, Y) 
        \qquad \text{ for all  }\, k 
        } 
      \ar[ru] \ar[rd] \\
      & j.
      }}
  $$
  Then, the classifying topos $\cat B G$ is equivalent to
  the category $\cat S$ of disjoint sums of objects of 
  $\cat B_0$. 

  The same result holds for progroupoids.
\end{lemma}

\begin{proof}
  This results from \cite[3.3.6]{leroy}, which states that 
  the projective $2$-limit of a
  \emph{strict} projective system of classifying toposes of groupoids
  (strict means that the inverse image functors are fully faithful:
  this is the case here)
  is the topos of sheaves on $\limind^{\dcat{Cat}} \cat B G_i$ for the
  canonical topology and from 
  \cite[2.4.2, 2.4.3]{leroy} which proves that 
  the category obtained from $\limind^{\dcat{Cat}} \cat B G_i$ by
  adding disjoint sums is a topos.
\end{proof}

\begin{defi}
  A \d{base point} of a topos $\T$ is a topos morphism 
  $\cat{Set} \lra \SLC\T$. 
\end{defi}
\begin{theo}
  Let $\T$ be a locally connected topos and $P$ a set of
  points of $\SLC\T$ (\ie, base points of $\T$), 
  at least one in each connected
  component of $\SLC\T$. 
  Let $G_R$ be the groupoid whose objects are elements of $P$ seen as
  points of $\LC(\T,R)$ and whose morphisms are topos $2$-isomorphisms
  $ \smash{
    \xymatrix @C=2ex { 
      *!!<0pt,0.7ex>+!R!<-1.2ex,0pt>
      {\cat{Set}}
      \ar@/^1ex/[r]|*{\vrule height 1.5pt depth 1.5pt width 0pt}="A"
      \ar@/_1ex/[r]|*{\vrule height 1.5pt depth 1.5pt width 0pt}="B"
      \ar@{=>}"A";"B"
      & 
      *!!<0pt,0.7ex>+!L!<+1.2ex,0pt>
      {\LC(\T,R)}
      } 
    }
  $
  and let 
  $\pi_1(\T, P) := \smash{\proobj _{R \in J(*)} G_R}$.
  Then $\pi_1(\T,P)$ is a strict progroupoid and 
  $\SLC\T \equiv \cat B \pi_1(\T,P)$.
\end{theo}
\begin{proof}
  To simplify, one may assume that $\T$ is connected. 
  But as the inverse image functors of the morphisms 
  $\T \lra \SLC\T$ or $\T \lra \LC(\T,R)$ 
  preserve any decomposition 
  $* = X_1 \amalg X_2$
  and are fully faithful, $\SLC\T$ and $\LC(\T,R)$ are also
  connected. 
  We may then assume that $P$ is a single point
  and that the $G_R$ are groups. We have 
  $\LC(\T,R) \equiv \cat B G_R$, thus 
  $\SLC\T \equiv \cat B \pi_1(\T,P)$.

  The projective system $(G_R)_{R \in J(*)}$ is clearly
  filtering. It remains to prove that the morphisms 
  $f: G_R \lra G_{R'}$ are surjective: this will come from the fact
  that 
  a group morphism
  $f: G' \lra G$ is onto 
  provided $f^*: \cat B G \lra \cat B G'$ is fully faithful.
  Indeed, as $G$ is a $G$-set, one has then a bijection 
  $G = \aut _{G} G \lra \aut_{G'} G = \aut _{\im f} G$.
  But an $\im f$-automorphism of $G$ is entirely determined by a
  permutation 
  $\sigma \in \SymetricGroup_{G/\im f}$ of the orbits and by an
  automorphism of each orbit; but as the automorphism group of each orbit
  is isomorphic to $\im f$, one has 
  $$ \aut_{G'} G = \aut_{\im f} G = 
  \Biggl( \prod_{G/\im f} \im f \Biggr)
  \psd 
  \SymetricGroup_{G/\im f}^{\mathstrut}.
  $$
  But as $G \lra \aut_{G'} G$ is a bijection, one easily sees that
  $\im f = G$. 
\end{proof}
\begin{rem}\label{rem:omissionofbasepoints}%
  By abuse of notation, we shall write 
  $\pi_1\T$ instead of 
  $\pi_1(\T, P)$ when $P$ is the set of all base points of $\T$. 
\end{rem}
\begin{rem}
  All this may be generalized to profinite groupoids in the following
  way. 
  Let $\T$ be a locally connected topos. 

  A \d{finite object} is an object of the form $I_\T$ for some finite
  set $I$, \ie, it is a disjoint sum of a finite number of copies of
  the final object. 
  A \d{locally constant finite object} is a locally constant object
  $F$ such that whenever $U$ trivializes $F$, the restriction 
  $F|_U$ is finite. 
  A \d{finite Galois object} is a locally constant finite Galois
  object. 
  A \d{finite Galois topos} is a topos generated by its finite Galois
  objects. 

  The full subcategory $\SLCF \T$ of disjoint sums of 
  locally constant finite objects of $\T$ is a finite Galois topos
  and the
  inclusion functor $\SLCF \T \lra \T$ is the inverse image functor of
  a topos morphism; this topos morphism is $2$-universal among those
  to finite Galois toposes. 
  If we let $\LCF(\T,R)$ denote the category of disjoint sums of
  locally constant finite objects of $\T$ trivialized by 
  the connected elements of a sieve $R$, then 
  $\SLCF\T = \dlimpro_{R \in J(*)} \LCF(\T,R)$.

  A \d{finite groupoid} is a groupoid whose hom-sets are finite. 
  A \d{profinite groupoid} is a filtered strict projective system of
  finite groupoids. The finite Galois toposes are the strict
  projective $2$-limits of classifying toposes of finite groupoids. 

  A \d{profinite base point} of $\T$ is a topos morphism 
  $\cat{Set} \lra \SLCF\T$. If $P$ is a set of profinite base points
  of $X$, at least one in each connected component of $\SLCF\T$, and
  if $G_R$ is the groupoid whose objects are elements of $P$ seen as
  points of $\LCF(\T,R)$ and whose morphisms are topos $2$-isomorphisms
  $ \smash{
    \xymatrix @C=2ex { 
      *!!<0pt,0.7ex>+!R!<-1.2ex,0pt>
      {\cat{Set}}
      \ar@/^1ex/[r]|*{\vrule height 1.5pt depth 1.5pt width 0pt}="A"
      \ar@/_1ex/[r]|*{\vrule height 1.5pt depth 1.5pt width 0pt}="B"
      \ar@{=>}"A";"B"
      & 
      *!!<0pt,0.7ex>+!L!<+1.2ex,0pt>
      {\LCF(\T,R)}
      } 
    }
  $,
  then 
  $\pi_1(\T, P) := \proobj _{R \in J(*)}$
  is a strict finite progroupoid and 
  $\SLCF\T \equiv \cat B \pi_1(\T,P)$.

  If $\cat C$ is a Galois category and $\f F : \cat C \lra \cat{set}$
  a fiber functor (see \cite{SGA1}), then the category of disjoint
  sums of objects of $\cat C$ is a connected finite Galois topos and
  $\f F$  induces a profinite base point, still denoted $\f F$, then 
  $\wh \pi_1( \T, \{\f F\})$ is the profinite fundamental group of
  $\cat C$.
  
  A \d{quotient groupoid} of a groupoid $G$ is a groupoid morphism $f:
  G \lra H$ such that $f$ induces a bijection between the sets of
  objects of $G$ and $H$ and such that the maps $\hom_G(x,y) \lra
  \hom(fx, fy)$ be surjective for all $x$, $y \in \ob G$. The
  \d{profinite completion} of a groupoid $G$ is the progroupoid 
  $\wh G$ of its finite quotients, the profinite completion of a
  progroupoid is the progroupoid of the finite quotients of its
  components. One may check that
  $\wh\pi_1\T$ is the profinite completion of $\pi_1\T$.
\end{rem}

\section{Fundamental group of a stack}

By ``stack'', we shall mean Deligne--Mumford stack
\cite{DM}, \cite{vistoli}, unless otherwise specified.
The following proposition enables us to apply the preceeding
theory to define the fundamental group of a stack. 

\begin{prop}
  The topos of étale sheaves on an algebraic stack is
  locally connected. 
\end{prop}
\begin{proof}
  The first lemma shows that it is generated by connected
  schemes, the second that these are connected as sheaves. 
\end{proof}
\begin{lemma}
  Let 
  $\fx F \dlra \fx G$ 
  be two distinct morphisms of étale sheaves on an algebraic stack
  $X$. Then, there exists a connected scheme $V$, an étale
  morphism $V \lra X$ (enabling us to see $V$ as a sheaf on
  $X$)
  and a
  morphism 
  $V \lra \fx F$ such that the two morphisms
  $V \lra \fx F \dlra \fx G$ 
  be distinct. 
\end{lemma}
\begin{proof}
  There exists an étale morphism $U \lra X$ such that the maps
  $\fx F U \iso \hom(U,\fx F) \dlra \hom(U,\fx G) \iso \fx G U$
  be distinct. If $(U_i)_{i\in I}$ are the connected components of
  $U$, then 
  the maps 
  $\prod \hom(U_i, \fx F) \dlra \prod \hom(U_i, \fx G)$ are distinct,
  hence there exists $i \in I$ such that the maps
  $\hom(U_i, \fx F) \dlra  \hom(U_i, \fx G)$ be distinct.
\end{proof}
\begin{lemma}
  Let $X$ be a stack and $U \lra X$ an étale morphism from a
  connected scheme. Then~$U$ is connected as an étale sheaf
  on $X$.  
\end{lemma}
\begin{proof}
  Let $U \lra X$ be an étale morphism from a connected scheme $U$, 
  seen as a presheaf 
  $T \lmt \hom_X(T,U)$ on $X\et$, 
  which is actually a sheaf.
  Assume $U = \fx F_1 \amalg \fx F_2$ and let us prove that either 
  $\fx F_1 = \nothing$ or $\fx F_2 = \nothing$.
  
  1. It suffices to prove that for any $Y \lra X$ in $X\et$ and for
  any geometric point 
  $x$ of $Y$, one has 
  $(\fx F)_x = \nothing$. Indeed, 
  $ 
  \{ 
  \,
  p^* : \cat{Sh}\, Y \lra \cat{Set}
  \text{ where }
  p
  \text{ is a geometric point of }
  Y
  \,
  \}
  $
  is a faithful family of fiber functors on $\cat{Sh}\,Y$, thus 
  \begin{align*}
    \{ 
    \,
    \xymatrix{ \cat{Sh}\,X \ar[r]^{f^*} & \cat{Sh}\, Y \ar[r]^{p^*} & \cat{Set} }
    \text{ where }
    &
    \xymatrix{ Y \ar[r]^f & X }
    \in X\et
    \text{ and }
    \\ &
    p
    \text{ is a geometric point of }
    Y
    \,
    \}    
  \end{align*}
  is a faithful family of fiber functors on $\cat{Sh}\,X$.

  2. If $x$ is a geometric point of $U$, then $U_x = *$ and 
  $U_x = (\fx F_1)_x \amalg (\fx F_2)_x$, thus 
  $(\fx F_1)_x = *$ and $(\fx F_2)_x = \nothing$ (or the converse).

  3. If $x : \spec k \lra U$ is a geometric point of $U$ such that 
  $(\fx F_1)_x = *$  and $(\fx F_2)_x = \nothing$, and if 
  $x' : \spec k' \lra U$ is a geometric point of $U$ above $x$, \ie, 
  $$ \xymatrix{\spec k' \ar[rd]^{x'} \ar[d] \\
    \spec k \ar[r]_x & U, } $$
  then there are maps 
  $U_{x'} \lra U_x$, 
  $(\fx F_1 )_{x'} \lra (\fx F_1)_x$ and 
  $(\fx F_2 )_{x'} \lra (\fx F_2)_x$, thus 
  $(\fx F_2 )_{x'} = \nothing$ and $(\fx F_1)_x=*$.
  As a result, the condition 
  ``$(\fx F_1)_x = *$ and $(\fx F_2)_x = \nothing$'' does not depend
  on the choice of a geometric point $x$ above a given point $x_0 \in
  X$. 

  4. The condition ``$(\fx F_1)_x = *$ and $(\fx F_2)_x = \nothing$''
  is open. Indeed, if 
  $(\fx F_1)_x \neq \nothing$, the, as
  $(\fx F_1)_x = \limind \fx F_1 (V)$, there exists 
  $V \lra U$ étale such that $\fx F_1(V)\neq \nothing$, and thus 
  $(\fx F_1)_y \neq \nothing$ for any $y$ in the image of 
  $V \lra U$ (which is open). 
  Similarily, the condition 
  ``$(\fx F_2)_x = *$ and $(\fx F_1)_x = \nothing$''
  is open, thus the condition 
  ``$(\fx F_1)_x = *$ and $(\fx F_2)_x = \nothing$''
  is closed. 
  As $U$ is connected, one has 
  $
  \forall x \in U 
  \quad 
  (\fx F_1)_x = * \text{ and } (\fx F_2)_x = \nothing
  $
  (or the converse).

  5. Let $V \lra X$ be an object of 
  $X\et$ and $x : \spec k \lra V$ be a geometric point of $V$. In the
  following diagram, the bold arrow are étale
  $$ \xymatrix{
    & U_x \ar[rr] \ar[ld] && 
    \smash[b]{U \Times_X V}
    \vrule width 0pt height 0pt depth 4pt
    \ar[ld]\ar[rd]^p \\
    \spec k \ar[rr]_x && V \ar[rd] && U \ar[ld] \\
    &&& X. }
  $$
  One may identify $U_x$ with the set of points $y$ of $U
  \times_X V$ over $x$.
  $$ \xymatrix{ 
    && \smash[b]{U \Times_X V} \ar[ld]\ar[rd] \\
    \spec k \ar@{.>}[rru]^y \ar[r]_x & V \ar[rd] && U \ar[ld] \\ 
    && X } $$
  As the decomposition 
  $U = \fx F_1 \amalg \fx F_2$ is preserved by pull-backs, 
  $U_x = \coprod_{y \in U_x} U_y$
  gives 
  $(\fx F_2)_x = \coprod_{y \in U_x} (\fx F_2) _y$. 
  But if 
  $y \in U_x$, one has a map 
  $(\fx F_2)_y \lra (\fx F_2)_{p(y)} = \nothing$, 
  thus 
  $(\fx F_2)_y = \nothing$.
\end{proof}

\begin{rem}
  This result also holds for the topos of étale sheaves on
  Artin stacks.
\end{rem}

\begin{rem}
  One may also define the profinite fundamental group of a
  stack $X$ using étale coverings, \ie, finite étale stack
  morphisms $Y \lra X$ instead of locally constant finite
  sheaves: see \cite{Z:ThVK}.
\end{rem}
\section{Comparison theorem}

We shall now compare the fundamental group of a stack over
$\CC$ and that of the associated (generalized) orbifold. 

\begin{defi}
  A \d{groupoid in a category} $\cat C$ (with finite projective
  limits) is the datum of 
  two objects $G_1$ and $G_0$, and of morphisms 
  $s, t: G_1 \lra G_0$, 
  $m : G_1 
  \times_{s, G_0, t} 
  G_1
  \lra G_1$, 
  $e: G_0 \lra G_1$ 
  and $i: G_1 \lra G_1$
  such that, for any object $T$,
  $(\hom(T,G_1), \hom(T,G_0), s, t, m, e, i)$ be a groupoid,
  where the maps $s$, $t$, $m$, $e$, $i$ are respectively 
  source, target, composition, identity and inverse. 
  A \d{topological groupoid} is a groupoid in the category of
  topological spaces.

  An \d{étale groupoid} is a topological groupoid whose
  source morphism (and hence, all of whose structure
  morphisms) is étale. 
  A \d{proper groupoid} is a topological groupoid such that
  the morphism
  $(s,t): G_1 \lra G_0 \times G_0$ be proper.
  One may also consider étale or proper groupoids in the category of
  schemes. 

  A \d{generalized orbifold} is an étale proper groupoid
  (see \cite{MP} for the notion of orbifold).
\end{defi}

\begin{prop}
  Let $X$ be a stack of finite type over $\CC$ 
  and $R \dlra U$ a presentation of $X$.
  The groupoid $R(\CC) \dlra U(\CC)$ is étale and proper. 
  We call it a (generalized) \d{orbifold associated to}~$X$.
\end{prop}
\begin{proof}
  The groupoid is clearly étale and properness results from 
  \cite[3.22]{LMB}.
\end{proof}

\begin{defi}
  An \d{equivariant sheaf} on an étale groupoid 
  $G_1 \dlra G_0$
  is a sheaf $\fx F$ on $G_0$ and a morphism 
  $\alpha : G_1 
  \mathbin{  {}_s{\mathop{\times}\limits_{G_0}} }
  \fx F \lra \fx F$
  such that 
  $$ \raisebox{.5\depth}{\xymatrix{
      G_1 \Times_{G_0} G_1 \Times_{G_0} \fx F \ar[r]^-{1\times\alpha}
      \ar[d]_-{m\times 1} & G_1 \Times_{G_0} \fx F \ar[d]^-\alpha \\
      G_1 \Times_{G_0}\fx F \ar[r]_\alpha & \fx F}}
  $$
  $$
  \raisebox{.5\depth}{\xymatrix{\relax
      \fx F = \smash{G_0 \Times_{G_0} \fx F}\vrule width 0pt depth 1ex
      \ar@{=}[rd]\ar[r]^-{e\times1} &
      G_1 \Times_{G_0} \fx F \ar[d]^\alpha \\ & \fx F}}
  \qquad
  \raisebox{.5\depth}{\xymatrix{ 
      G_1 \Times_{G_0} \fx F \ar[r]^\alpha \ar[d] & \fx F
      \ar[d] \\ G_1 \ar[r]_b & G_0. }} $$ 
  Equivalently, one may require that 
  $(\pr_2, \alpha) : G_1 \times G_0 \fx F \dlra \fx F$ 
  be the source and target maps of a groupoid such that 
  $$ \xymatrix{
    G_1 \times \fx F \ar[d] \dar[r] & \fx F \ar[d] \\ 
    G_1 \dar[r] & G_0 } $$
  be a groupoid morphism. 

  The equivariant sheaves and their morphisms form a topos,
  denoted $\cat{Sh}(G_1 \dlra G_0)$.  

  The \d{fundamental groupoid} of an étale groupoid 
  $G_1 \dlra G_0$
  is defined to be 
  $$ 
  \pi_1 ( G_1 \dlra G_0 ) := 
  \pi_1 \cat{Sh}(G_1 \dlra G_0).
  $$
\end{defi}

\begin{rem}
  \cite{vistoli} shows that if $R\dlra U$ is a presentation
  of a stack $X$, then $\cat{Sh} \, X \equiv \cat{Sh} (R
  \dlra U)$.
\end{rem}
\begin{rem}
  As there is an equivalence of categories between the
  category of coverings of an orbifold and that of étale
  coverings of one of its presentations, this definition
  generalizes that of the fundamental group of an orbifold
  \cite[13.2.5]{thurston}.
\end{rem}

\begin{theo}
  Let $X$ be a separated algebraic stack, locally of finite
  type over $\CC$, such that there exists a surjective étale
  and \emph{finite} morphism $U \lra X$ from a scheme.  We
  let $R \dlra U$ be the groupoid defined by this morphism.
  Then one has
  $$\wh\pi_1 X
  \equiv 
  \wh\pi_1 \bigl( R(\CC)\dlra U(\CC) \bigr).$$
\end{theo}

\begin{proof}
  \cbstart
  We shall identify locally constant finite
  étale sheaves on the groupoid $R \dlra U$ with its equivariant étale
  coverings. 
  It suffices to show that the functor 
  \begin{equation}
    \label{eq:functcompgr}
    \definefunction
    {\SLCF\cat{Sh}(R\dlra U)}%
    {\SLCF\cat{Sh}\bigl(R(\CC)\dlra U(\CC)\bigr)}%
    {\XYMATRIX{
        R \Times_U F \ar[d]
        \ar@<+2pt>[r]^-{\pr_2}
        \ar@<-2pt>[r]_-{\alpha}
        & F \ar[d]^f \\
        R
        \ar@<+2pt>[r]^s
        \ar@<-2pt>[r]_t
        & U
        }}%
    {\XYMATRIX{
        R(\CC) \Times_{U(\CC)} F(\CC) \ar[d]
        \ar@<+2pt>[r]
        \ar@<-2pt>[r]
        & F(\CC) \ar[d]^f \\
        R(\CC)
        \ar@<+2pt>[r]
        \ar@<-2pt>[r]
        & U(\CC)
        }}%
  \end{equation}
  is an equivalence of categories. 

  Let $(F, f, \alpha)$ and $(F', f', \alpha')$ be two étale coverings
  of $R \dlra U$. Let us first show that the map 
  \begin{equation}
    \label{eq:mapproofcompth}
    \raisebox{.5\depth}{\xymatrix{
        \hom\bigl(
        (F,f,\alpha), 
        (F', f', \alpha')
        \bigr)
        \ar[d]
        \\
        \hom\bigl(
        (F(\CC),f(\CC),\alpha(\CC)), 
        (F'(\CC), f'(\CC), \alpha'(\CC))
        \bigr)
        }}
  \end{equation}
  is onto. Indeed, any $\phi_0: F(\CC) \lra F'(\CC)$ with 
  $ f(\CC) = f'(\CC) \circ \phi_0$
  is of the form 
  $\phi_0 = \phi(\CC)$ for some 
  $\phi : F \lra F'$ with 
  $f = f ' \circ \phi$
  (see \cite[XII 5.1]{SGA1}); 
  further, as 
  the diagram 
  $$ \xymatrix{ 
    R(\CC) \Times_{U(\CC)} F(\CC) 
    \ar[r]^-\alpha \ar[d]_{1\times\phi} & 
    F(\CC) \ar[d]^\phi \\
    R(\CC)\Times_{U(\CC)} F'(\CC) 
    \ar[r]_-{\alpha'} & F'(\CC) } $$
  commutes and as that the functor 
  \begin{equation}
    \label{eq:functproofcompesptopassoc}
    \definefunction
    {\cat{Sch}/\CC}{\cat{Top}}{X}{X(\CC)}
  \end{equation}
  is faithful, we see that the diagram 
  $$ \xymatrix{ 
    R \Times_U F 
    \ar[r]^\alpha \ar[d]_{1\times\phi} 
    & F \ar[d]^\phi \\
    R\Times_U F' \ar[r]_{\alpha'} & F' } $$
  commutes.
  The faithfulness of \eqref{eq:functproofcompesptopassoc}
  also implies that the map \eqref{eq:mapproofcompth} is
  one-to-one. Hence, the functor \eqref{eq:functcompgr} is fully
  faithful. 

  The functor \eqref{eq:functcompgr} is essentially surjective. 
  Indeed, if 
  $$ \xymatrix{ 
    R(\CC) \Times_{U(\CC)} F(\CC) \ar[d]_{1\times f_0}
    \ar@<+2pt>[r]^-{\pr_2}
    \ar@<-2pt>[r]_-{\alpha_0}
    & F(\CC) \ar[d]^{f_0} \\
    R(\CC)
    \ar@<+2pt>[r]^s
    \ar@<-2pt>[r]_t
    & U(\CC)
    }
  $$
  is an equivariant étale cover of $R(\CC)\dlra U(\CC)$, then 
  by \cite[XII 5.1]{SGA1}
  there exists an étale cover 
  $f : F \lra U$ such that 
  $F(C) \iso F_0$ (we shall thus identify $F(\CC)$ and $F_0$) and 
  $f(\CC) = f_0$. 
  As 
  $U \lra X$ is finite, so is 
  $(R \times_U F)(\CC) \iso R(\CC) \times_{U(\CC)} F_0 \lra U(\CC)$
  and thus 
  by \cite[XII 5.1]{SGA1}
  there exists 
  $\alpha : R \times_U F \lra F$ such that $\alpha(\CC) = \alpha_0$. 
  One may then show as above that 
  $(F, f, \alpha)$ is an equivariant covering, 
  again from the fact that $(F_0, f_0, \alpha_0)$ is one and the fact that 
  the functor \eqref{eq:functproofcompesptopassoc} is faithful. 
  \cbend
\end{proof}

\section{Simplicial schemes}

We shall now show that our notion of fundamental groupoid of
a stack coincides with that of \cite{oda}.

\begin{defi}
  Let $\Delta$ be the \d{category of simplices}, \ie, the
  category whose objects are the sets $[n] = \{0,\dots,n\}$
  and whose morphisms are nondecreasing maps, and 
  $\Delta_n$ the full subcategory of $\Delta$ containing 
  $[0],\dots,[n]$.

  Let $\T$ be a locally connected topos and denote by $\pi$
  its ``connected components'' functor. 

  A \d{simplicial object} in $\T$ is a functor 
  $\Delta\op \lra \T$; we shall denote 
  $\Delta\op\T$ the category of simplicial objects of $\T$. 
  The restriction functor 
  $\Delta\op\T \lra \Delta_n\op\T$ has left and right
  adjoints, called $n$th \d{skeleton} and \d{coskeleton}
  denoted $\sk_n$ and $\cosk_n$.

  A \d{hypercovering} is a simplicial object $X_\bullet$
  such that the morphisms 
  \begin{align*}
    X_0       &\lra e \\
    X_{n+1}   &\lra \bigl(\cosk_{n}(X_\bullet)\bigr)_{n+1}
  \end{align*}
  be epic. 

  If $I$ is a set and $X$ an object of $\T$, set 
  $I \otimes X := \coprod_{i\in I} X$.
  If $I_\bullet$ is a simplicial set ans $X_\bullet$ a
  simplicial object, set
  $I_\bullet \otimes X_\bullet := (n \lmt I_n \otimes X_n).$

  A \d{strict homotopy} between two morphisms 
  $f,\,g : X_\bullet \lra Y_\bullet$ is a morphism 
  $[1] \otimes X_\bullet \lra Y_\bullet$ whose restrictions
  to 
  $[0]\otimes X_\bullet = X_\bullet$ are $f$ and $g$. 
  We call \d{homotopy} the equivalence relation 
  generated by strict homotopy. 
  We denote  $\cat{HR}\T$ the category of hypercoverings of
  $\T$ and their morphisms up to homotopy. One may show that
  $(\cat{HR}\T)\op$ is filtering \cite{AM}.

  Finally, let $\H$ denote the category of CW-complexes and
  their morphisms modulo homotopy (which may be seen as a
  category whose objects are simplicial sets) and 
  $\pro\H$ the category of its proobjects. 
\end{defi}

\begin{defi}
  The \d{homotopy type} of a locally connected topos $\T$ is
  $$ \{\T\} := \proobj_{Y_\bullet \in \cat{HR}\T}
  \pi Y_\bullet 
  \in \pro\H.
  $$
  Its \d{fundamental groupoid} is 
  $$ \pi_1 \{\T\} := \proobj_{Y_\bullet \in \cat{HR}\T}
  \pi_1 \pi Y_\bullet
  $$
  (with the same convention as in remark
  \ref{rem:omissionofbasepoints}
  for the omission of base
  points). 
\end{defi}

\begin{prop}
  If $\T$ is a locally connected topos, then 
  $$ \cat B \pi_1 \{ \T \} \equiv \SLC \T. $$
\end{prop}
\begin{proof}
  One has
  \begin{align*}
    \cat B \pi_1 \{ \T \} 
    &= \cat B \pi_1 \proobj_{U_\bullet \in \cat{HR}\,\T} 
    \pi U_\bullet
    &&\text{by the definition of $\{\T\}$}
    \\
    &= \cat B \proobj_{U_\bullet \in \cat{HR}\,\T} 
    \pi_1 \pi U_\bullet 
    &&\vtop{\hbox{by the definition of the fundamental}%
      \vskip -4pt
      \hbox{group of a pro-simplicial set}}
    \\
    &= \dlimpro_{U_\bullet \in \cat{HR}\,\T} 
    \cat B \pi_1 \pi U_\bullet 
    &&\vtop{\hbox{by the definition of the 
        classifying}%
      \vskip -4pt
      \hbox{space of a progroup}}
    \\
    &\equiv \dlimpro_{U_\bullet \in \cat{HR}\,\T} \LC(\T, U_0) 
    &&\text{from \cite[10.6]{AM}} \\
    &\equiv \dlimpro_{R \in J(*)} \LC(\T, R). 
  \end{align*}
  To prove the last equivalence, it suffices to show that
  any sieve $R$ is finer than a sieve generated by a single
  object $U_0$. One may be tempted to set $U_0 =
  \coprod_{(U\to*)\in R} U$, but this coproduct need not be
  small. Instead, as there exists a \emph{set} $\cat K
  \subset \ob\T$, generating $\T$, one may set
  $$ U_0 = \coprod_{\substack{
      C \in \cat K \text{ such that} \\
      C\to* \text{ factors} \\
      \text{through a } (U\to*)\in R
      }}
  C.
  $$
\end{proof}

\section{Short exact sequences}

\begin{theo}
  Let $\T$ be a connected locally connected topos, $x$ a
  point of $\SLC\T$ and $\f F : \SLC\T \lra \cat{Set}$ the
  corresponding fiber functor.  Let~$Y$ be a Galois object
  of $\SLC\T$.
  \begin{theoenum}
  \item The category $\SLC\T /Y$ is a Galois topos and a
    point $a \in \f F(Y)$ defines a fiber functor $\SLC\T
    \lra \cat{Set}$.
  \item The functor 
    $$\definefunction{\SLC\T }{(\SLC\T )/Y}%
    X{(X\times Y\lra Y)} 
    $$
    induces a topos morphism
    $(\SLC\T )/Y \lra \SLC\T $.
  \item There is a short exact sequence
    $$ 1 \lra \pi_1(Y, a) \lra \pi_1(\T, x) 
    \lra \aut Y \lra 1. $$
  \end{theoenum}  
\end{theo}
\begin{proof}
  One may assume that $\T$ is the classifying topos of a
  progroup and that $\f F$ is the forgetful functor
  $G\text{-}\cat{Set} \lra \cat{Set}$.  The following lemma
  identifies $G\text{-}\cat{Set}/Y$ as the classifying topos
  of a subprogroupoid $H$ of~$G$, thereby proving $(i)$ and
  $(ii)$.
  \begin{lemma}
    Let $G$ be a progroup, $Y$ a connected $G$-set, 
    $a \in Y$ and $H$ the stibilizer of $a$. 
    Then, there is an equivalence of categories 
    $$ G\text{-}\cat{Set} / Y \equiv H\text{-}\cat{Set}. $$
  \end{lemma}
  \begin{proof}
    Consider the functor 
    $$\definefunction
    {G\text{-}\cat{Set}/Y}
    {H\text{-}\cat{Set}}
    {(\xymatrix{Z \ar[r]^\pi & Y})}
    {\pi^{-1}(\{a\}).}
    $$

    Let $Z$ and $Z'$ be objects of 
    $G\text{-}\cat{Set}/Y$. The map
    $$ \hom_{Y,G}(Z,Z') \lra \hom_H\Bigl(
    \pi^{-1}\bigl(\{a\}\bigr), 
    \pi'{}^{-1}\bigl(\{a\}\bigr)
    \Bigl) $$
    is injective, because a morphism
    $f : Z \lra Z'$ is entirely determined by its fiber
    $f_a$. Indeed, as $Y$ is connected, for any $z \in Z$,
    there exists $g \in G$ such that $\pi(z) = g \cdot a$, 
    thus $f(z) = f(gg^{-1}\cdot z)  = 
    g \cdot f(g^{-1}\cdot z))$, for 
    $g^{-1} \cdot z \in \pi^{-1}(\{a\})$.
    
    To show that the map is onto, let 
    $f_a : \pi^{-1}(\{a\}) \lra \pi'{}^{-1}(\{a\})$
    be an $H$-morphisms and set 
    $f(z) = g \cdot f_a( \gamma^{-1} \cdot z )$ 
    if 
    $\pi(z) = g \cdot a$. 
    This is well-defined, for if 
    $\pi(z) = g \cdot a = g' \cdot a$, 
    then 
    $g^{-1} g' a =  a$
    thus  $h = g^{-1} g' \in H$, 
    but as $f_a$ is $H$-equivariant, 
    $h^{-1} f_a (g^{-1} \cdot z ) =
    f_a( h^{-1} g^{-1} z )$, 
    thus 
    $g f_a( g^{-1} z ) = g' f_a( g'{}^{-1} z )$.
    Further, the map $f$ is equivariant: 
    $g' \cdot f(z) =  
    g' g \cdot f_a( (g' g)^{-1} \cdot g' z ) =
    f(g' \cdot z)$.

    Our functor is thus fully faithful. As for essential
    surjectivity, let $U$ be an $H$-set and define a
    $G$-action on 
    $Y \times U$ as follows: 
    for all $y \in Y$, choose 
    $\gamma_y \in G$ such that 
    $y = \gamma _y \cdot a$ (one may set $\gamma_a=1$) and
    set 
    $g \cdot(y,u) = (g \cdot y, 
    \gamma_{g y} ^{-1} g \gamma_y u)$. 
    This is clearly a $G$-action and the action
    of $H$ on the fiber is the desired one. 
  \end{proof}

  The following lemma shows that 
  $\pi_1(Y, a)$ is a normal subgroup of 
  $\pi_1(\T, x)$ and identifies the quotient, thereby
  establishing 
  $(iii)$.
\end{proof}

\begin{lemma}
  Let $G$ be a progroup, 
  $Y$ a Galois $G$-set, 
  $a \in Y$ and $H$ the stabilizer of $a$. 
  Then, the subprogroup $H$ is normal and there is an
  isomorphism
  $$ G / H \iso \aut Y.$$
\end{lemma}
\begin{proof}
  Let us first assume that $G$ is a group. 
  As $Y$ is Galois, the map 
  $$ \Psi : \definefunction 
  {Y \times \aut Y}
  {Y \times Y}
  {(y,\phi)}
  {(y,\phi(y))}
  $$
  is bijective so we may define 
  $$ \Phi : \definefunction
  {G}
  {\aut Y}
  {g}
  {\phi \text{ such that } \Psi(a,\phi) = (a, g\cdot a).}
  $$
  This map is onto: as $Y$ is connected, an automorphism 
  $\phi \in \aut Y$ is entirely determined by 
  $\phi(a)$, but there exists a 
  $g \in G$ such that 
  $\phi(a) = g \cdot a$, thus 
  $\phi = \Phi(g)$. This map is clearly a group morphism
  and its kernel is readily seen to be $H$. 
  
  In the case $G = \proobj G_i$ is a progroup, one has
  short exact sequences 
  $$ 1 \lra H_i \lra G_i \lra \aut Y \lra 1$$
  thus, by \cite[I.8.9.2]{SGA4}, 
  $$ 1 \lra H \lra G \lra \aut Y \lra 1.$$
\end{proof}

\begin{coro}
  Let $X$ be a connected scheme, $G$ a finite group actiong on $X$ and
  $x$ a geomertic point of $X$. Then, there is a short exact sequence 
  $$ 1 \lra \wh\pi_1 (X, x) \lra \wh\pi_1( [X / G], x) \lra G \lra 1. $$
\end{coro}
\begin{proof}
  Let us apply the theorem to
  $\T = \SLCF \cat{Sh}(G \times X \dlra X)$, 
  the topos
  of disjoint sums of $G$-equivariant étale coverings of $X$ and the
  object $Y = G \times X$, where the action of $G$ and the morphism 
  $Y \lra X$ are 
  $$ \definefunction
  {G \times Y}
  {Y}
  {(h,g,x)}
  {(gh^{-1}, hx)}
  \qquad\qquad
  \pi : \definefunction
  {Y}{X}
  {(g,x)}
  {g \cdot x}
  $$

  To prove that 
  $G = \aut Y$, remark that the one-to-one morphism 
  $$ \Phi : \definefunction
  {G}
  {\aut Y}
  {h}
  { \Phi(h) : \definefunction{Y}{Y}
    {(g,x)}{(h^{-1}g, h \cdot x)}
    }
  $$
  is onto: indeed, if 
  $\phi \in \aut Y$, $x \in X$, 
  as 
  $\phi$ is an automorphism over~$X$, one may write 
  $\phi(1,x) = (h^{-1}, h \cdot x)$ with 
  $h \in G$: thus 
  $\phi$ and $\Phi(h)$ coincide at a point $(1,x)$ of $Y$, hence on
  the connected component 
  $\{1\}\times X$, hence on its orbit $Y$.

  To show that $Y$ is Galois, it suffices to remark that 
  \begin{align*}
    Y \Times_X (\aut X) &=
    \{ (g,h,x,y) \in G\times X\times G\times X \ :\ 
    g\cdot x = y \}
    \\
    Y \Times_X Y &=
    \{ (g,h,x,y) \in G\times X\times G\times X \ :\ 
    g\cdot x = h\cdot y \}.
  \end{align*}

  Finally, as 
  $\T/Y \equiv \SLCF\cat{Sh}\,X$, the short exact sequence of the
  theorem 
  $$ 1 \lra \wh\pi_1 \T/Y \lra 
  \wh\pi_1\T \lra
  G \lra 1 $$
  reads 
  $$ 1 \lra \wh\pi_1 X \lra 
  \wh\pi_1 [X/G] \lra
  G \lra 1. 
  $$
\end{proof}
\begin{rem}
  The same result holds for stacks quotiented by a finite group.
\end{rem}
\begin{coro}
  Let $X$ be a connected quasi-compact algebraic stack on a field $k$,
  let $\bar k$ be a separable closure of $k$ and let 
  $x$ be a point of $\bar X = X \tens_k \bar k$. 
  Then, one has a short exact sequence of profinite
  fundamental groups 
  $$ 1 \lra \wh\pi_1(\bar X, x) \lra 
  \wh\pi_1 (X,x) \lra 
  \gal(\bar k : k ) \lra 1. $$
\end{coro}
\begin{proof}
  By lemma \ref{lemma:dsfsdkfgkjjd1},
  for each Galois extension 
  $k':k$, 
  $X \tens_k k'$ is a Galois object in 
  $\SLCF\cat{Sh}\,X$
  and the short exact sequence of the theorem reads 
  $$ 
  1 \lra \wh \pi_1(X \tens_k k') \lra \wh\pi_1(X) 
  \lra \gal(k':k) \lra 1,
  $$
  thus, from \cite[I.8.9.2]{SGA4}, 
  $$ 
  1 \lra \limpro_{k':k} \wh \pi_1(X \tens_k k') 
  \lra \limpro_{k':k} \wh\pi_1(X) 
  \lra \limpro_{k':k} \gal(k':k) \lra 1,
  $$
  but lemma \ref{lemma:dsfsdkfgkjjd2} identifies the first term,
  yielding 
  $$ 1 \lra \wh\pi_1(X \tens_k k') 
  \lra \wh \pi_1(X) \lra \gal(\bar k : k) \lra 1$$
  as desired.
\end{proof}
\begin{lemma}\label{lemma:dsfsdkfgkjjd1}%
  Let $X$ be a connected algebraic stack over a field $k$. Then 
  $\smash{X \Otimes_k k' \lra X}$ is Galois. 
  and 
  $ \gal (k' : k ) \iso \aut_X (X \otimes_k k') $
\end{lemma}
\begin{proof}
  Let $s,t : R \dlra U$ be a presentation of $X$, 
  $\spec A_i$ affine open sets covering $U$, 
  $\spec B_{ijk}$ affine open sets covering 
  $s^{-1} \spec A_i \inter t^{-1} \spec A_j$.
  An automorphism of 
  $X \tens _k k'$ over $X$ induces 
  automorphisms of 
  $\spec k' \tens_k B_{ijk}$ over 
  $\spec B_{ijk}$
  and of 
  $\spec k' \tens_k A_i$ over $\spec A_i$, 
  compatible with $s$ and $t$; 
  each of these uniquely determines an automorphism of $k'$ over~$k$,
  but as the groupoid 
  $R \dlra U$ is connected, it is always the same automorphism of
  $k':k$. 

  To show that $X'$ is Galois, just compute:
  \begin{align*}
    X' \Times_X X' &=
    \left(
      X \Times_{\spec k} \spec k'
    \right)
    \Times_X
    \left(
      X \Times_{\spec k} \spec k'
    \right)
    \\
    &=
    X
    \Times_{\spec k}
    \left(
      \spec k' \Times_{\spec k} \spec k'
    \right)
    \\
    &=
    X
    \Times_{\spec k}
    \left(
      \spec k' \times \gal(k':k)
    \right)
    &&\text{for $k':k$ is Galois}
    \\
    &=
    X
    \Times_{\spec k}
    \aut_X (X \tens_k k').    
  \end{align*}
\end{proof}
\begin{lemma}\label{lemma:dsfsdkfgkjjd2}%
  Let $X$ be a connected
  algebraic stack over a field $k$, with a quasi-compact presentation,
  and let $\bar k$ be a separable closure of $k$. Then, there is an
  isomorphism of profinite groups 
  $$ \wh\pi_1( X \otimes_k \bar k) \iso \proobj_{\substack{
      k':k \text{ finite}\\\text{Galois extension}}}
  \wh\pi_1( X \otimes_k k' ). $$
\end{lemma}
\begin{proof}
  It suffices to prove that any étale covering of 
  $\bar X = X \tens_k \bar k$ comes from a covering of 
  $X' = X \tens_k k'$ for some \emph{finite} Galois extension 
  $k':k$. Let $s,t : R \dlra U$ be a presentation of~$X$, with $R$ and
  $U$ quasi-compact, 
  let 
  $\alpha: R\times_{\bar U} F \lra F$
  be an equivariant étale cover of 
  $\bar R \dlra \bar U$,
  let 
  $\spec A_i$ be a finite number of open affine sets covering~$U$, 
  let $\spec B_{ij}$ be a finite number of open affine sets covering 
  $\alpha^{-1} \spec A_i$. 
  The morphisms 
  $F \lra \bar U$ and $\alpha$ may be written, locally, 
  \begin{align*}
    \spec \dfrac{A_i \tens \bar k [ x_1, \dots, x_{n_i} ]}%
    {(f_{i,1}, \dots, f_{i,n_i})} 
    &\lra 
    \spec A_i\tens \bar k,
    \\
    \spec B_{ij} \Otimes_{A_i} 
    \dfrac{A_i \tens \bar k [ x_1, \dots, x_{n_i} ]}%
    {(f_{i,1},\dots,f_{i,n_i})} 
    &\lra
    \spec     
    \dfrac{A_i \tens \bar k [ x_1, \dots, x_{n_i} ]}%
    {(f_{i,1},\dots,f_{i,n_i})}
  \end{align*}
  and only involve a finite number of elements of 
  $\bar k$: thus, $F$ is defined over the finite Galois extension
  generated by these elements. 
\end{proof}
\bibliographystyle{amsplain}
\nocite{*}
\bibliography{1}
\end{document}